\documentclass{articles_perso}

\title{Instability of the magnetohydrodynamics system at small but finite Reynolds number}
\date{}

\begin{document}
\begin{abstract}
	The aim of this paper is to give a result concerning the stability properties of the solutions of magnetohydrodynamics equations at small but finite Reynolds numbers. These solutions are found using the \(\alpha\)-effect: this method gives us solutions which are highly oscillating spatially on the scale of the underlying flow but are growing on a larger scale depending on a parameter \(\varepsilon\). We prove instability results for a dense subset of initial velocity field of the flow at given Reynolds number.
	\end{abstract}
\begin{keywords}
	magnetohydrodynamics, stability, \(\alpha\)-effect
	\end{keywords}
\section{Introduction}
Dynamo theory first interested physicists almost one century ago, when Sir Joseph Larmor \cite{larmor1919could} wondered how a magnetic field could be generated by a rotating body like the Sun or the Earth. Since then, many results have come, both in mathematics and physics, to understand this kind of phenomena.

The physical model behind them is a coupled system derived from the interaction of a conductive fluid, which follows Navier-Stokes equations, and an electromagnetic field, which follows Maxwell equations and Ohm's law. The interaction of the two physic quantities comes from the fact that a magnetic field has an effect on a moving conductive fluid, which creates in return a magnetic field. This kind of system is known as magnetohydrodynamics (MHD).

There are mainly two outcomes for MHD systems: the magnetic energy can either collapse (exponentially), or stay/grow. In the latter case we say we have a dynamo, and it is the effect we observe in bodies like the Sun. However, the understanding of dynamo mechanisms is not simple. Indeed, the most ``intuitive'' way we imagine the conductive fluid in rotation and the magnetic field (that is a 2D rotating fluid in the kernel, and a magnetic field that looks like the one of a static magnet) is actually a well known anti-dynamo situation. Anti-dynamo theorems have been widely studied, and some summary of the knowledge in the domain has been made in \cite{gilbert2003dt}.

The kind of situation we want to study there is that of the second kind, that is dynamo situation (and more precisely growth). Positive results began to emerge only later in the 50's, even though it was already well known that it is what occured in the Sun. In the 90's M.M. Vishik \cite{vishik1989mfg} proved that fast dynamos (i.e. exponential growth) require Lagrangian chaos in the fluid. Those results spread the idea that the velocity field has to be chaotic enough to expect a dynamo effect. Because of this ``negative'' result, different methods were tried to find instability results, linked to singularities: Y.B. Ponomarenko \cite{ponomarenko1973theory} studied concentration phenomena to get dynamo result in the case where Laplace force is neglected (i.e. imposed velocity field) and his results were extended by Gilbert \cite{gilbert1988fast} and successfully experimentaly checked \cite{gailitis2004riga}. More realistic case of the Ponomarenko dynamo was studied by D. G\'erard-Varet and F. Rousset \cite{gerard2007shear}. 
All those results give only asymptotic results.

In this paper we will focus on another singular dynamo mechanism, called \emph{alpha effect}, first introduced by E.N. Parker \cite{parker1955hydromagnetic}, which is based on scale separation: we formally decompose the magnetic and the velocity fields into two parts: a fluctuating part, evolving on a turbulent scale, and a mean part, evolving on a larger scale. The mean value is broadly the spatial average over a sphere of radius much larger than the typical distance on which the field is turbulent, but much smaller thant the distance on which we want to observe the dynamo effect (in what follows, these distances will be \(1\) and \(1/\varepsilon\), so scale separation will take place at small \(\varepsilon\)). Basically, this gives a decomposition of the velocity and magnetic fields under the form
\begin{equation}
  \begin{split}
	u &= \tilde u(x,x/\varepsilon,t) + \bar u(x,t), \\
	b &= \tilde b(x,x/\varepsilon,t) + \bar b(x,t), \\
        \end{split}
	\end{equation}
where \(\tilde u\) and \(\tilde b\) have null mean value on the second variable. The idea behind \emph{alpha effect} is that the averaged induction term \(\overline{\nabla\wedge(\tilde u\wedge\tilde b})\) (we will specify those notations right below) that will appear in the MHD equations will permit us to obtain a destabilizing effect on the magnetic field.

The \emph{alpha effect} was first formally studied by G.O. Roberts \cite{roberts1972dynamo,roberts1970spatially}, in the framework of periodic flows. His results were discussed by H.K. Moffatt \cite{moffatt1978fge}, S. Childress \cite{childress1979alpha} and A.M. Soward \cite{soward1987fast} that considered a larger class of velocity fields and range of parameters.
This effect has been since experimentaly confirmed by R. Stieglitz and U. M\"uller \cite{stieglitz2001experimental}.

The aim of the present paper is to provide a rigorous justification of this mechanism, in a nonlinear framework. We will improve substantially former results of D. G\'erard-Varet \cite{gerardvaret2006osi,gerardvaret2007gafd}.

Before stating more precisely our results, we will first recall the system we will work with:

The generation of magnetic field \(B\) by an electrically conducting fluid of velocity \(U\) and magnetic diffusivity \(\nu_m\) follows the equations
\begin{equation}
	\left\{\begin{array}{l}
		\partial_tB - \nabla\wedge(U\wedge B) - \nu_m\Delta B = 0,\\
		\nabla\cdot B = 0.
		\end{array}\right.
	\end{equation}
These equations are derived from the Maxwell's equation and Ohm's law (see second section of \cite{fauve2003dynamo}, or \cite{gilbert2003dt} for more detailed derivation). Since other relevant physics quantities that appear in the Maxwell's equations can be obtained from \(B\), it is enough to work with only that.

Furthermore, the incompressible Navier-Stokes equations that drive the evolution of the conducting fluid under the force created by the magnetic field write
\begin{equation}
	\left\{\begin{array}{l}
		\rho\partial_tU + \rho U\cdot\nabla U = -\nabla p + \eta\Delta U + \frac1{\mu_0}(\nabla\wedge B)\wedge B + f,\\
		\nabla\cdot u = 0,
		\end{array}\right.
	\end{equation}
where \(p\) is the pressure field, \(f\) is an additional forcing term, \(\rho\) is the fluid mass density and \(\eta\) is the fluid shear viscosity (see \cite{fauve2003dynamo} again, keeping in mind the incompressible condition). The additional forcing term will be chosen to generate some steady velocity field (and thus will be independent of time), see below.

We study there oscillating field \(U\) of amplitude \(V\), oscillating on a space-scale \(L\), and with magnetic diffusivity \(\nu_m\). After proper nondimensionalization of time, space, pressure and magnetic field, we then have the full incompressible MHD system:
\begin{equation}
	\label{MHD}
	\left\{\begin{array}{l}
		\partial_tB - \nabla\wedge(u\wedge B) - \frac1{R_m}\Delta B = 0,\\
		\partial_tU +U\cdot\nabla U + \nabla p - \frac1{R_e}\Delta U = (\nabla\wedge B)\wedge B +f,\\
		\nabla\cdot B = 0,\\
		\nabla\cdot U = 0,\\
		\end{array}\right.
	\end{equation}
where \(R_m = \frac{L V}{\nu_m}\) is the magnetic Reynolds number, and \(R_e = \frac{\rho LV}{\eta}\) is the hydrodynamic Reynolds number.

Here we want to find a growing solution on a larger scale \(L/\varepsilon\) (where \(\varepsilon>0\)), using the \emph{alpha effect} introduced earlier. We will search this solution as a perturbation \(\begin{pmatrix} u \\ b\end{pmatrix}\) of a solution \(\begin{pmatrix} U_s \\ 0\end{pmatrix}\) of this system (and the forcing term \(f\) is chosen accordingly so that \(U_s\) is a solution of the Navier-Stokes equation):
\begin{equation}
	\left\{\begin{array}{l}
		U = u + U_s, \\
		B = b + 0.
		\end{array}\right.
	\end{equation}
and the system in terms of \(b\), \(u\) and \(U_s\) becomes
\begin{equation}
	\label{MHDPert}
	\left\{\begin{array}{l}
		\partial_tu + u\cdot\nabla U_s + U_s\cdot\nabla u - \frac1{R_e}\Delta u + \nabla p = (\nabla\wedge b)\wedge b - u\cdot\nabla u, \\
		\partial_tb - \nabla\wedge(U_s\wedge b) - \nabla\wedge(u\wedge b) - \frac1{R_m}\Delta b = 0,\\
                \nabla\cdot b = 0,\\
                \nabla\cdot u = 0.
		\end{array}\right.
	\end{equation}


The space variables for this system is denoted \(\theta\), and corresponds to the fluctuation scale \(L\), \(1\) after non-dimensionalization. In particular, we assume that
\begin{equation}
  U_s = U_s(\theta)\in\mc C^\infty(\mb T^3)
  \end{equation}
and has zero-mean.

Later on, we shall introduce a larger scale variable \(x = \varepsilon\theta\) corresponding to the large scale \(L/\varepsilon\) of the growing field.

Our main result is an instability theorem in the Liapunov sense for system \eqref{MHDPert}. To state our result, we introduce the set \(\mc H_0^\infty(\mb T^3)\) with his usual distance:
\begin{equation}
	d(U,U') = \sum_{k=0}^\infty\frac{\max\left(1,\|U-U'\|_{\mc H^k}\right)}{2^k},
	\end{equation}
which makes it a Fr\'echet-space. 

For \(T = (T_1,T_2,T_3)\in(\mb R_+^*)^3\), we denote
\begin{equation}
  \mc H_T^s = \mc H^s(\mb R/T_1\mb Z\times\mb R/T_2\mb Z\times\mb R/T_3\mb Z)
\end{equation}
the Sobolev space of order \(s\) on the torus
\begin{equation}
  \mb R/T_1\mb Z\times\mb R/T_2\mb Z\times\mb R/T_3\mb Z.
\end{equation}

We have the following theorem:
\begin{theoreme}
  Let \(s>5/2\) be a real number. There exists a dense subset \(\ms P\) of the unit ball of \(\mc H_0^\infty(\mb T^3)\) and \(R_m^0>0\) such that for all \(U_s\in\ms P\)  and \(R_m<R_m^0\), for all \(R_e>0\), 
the solution \(\begin{pmatrix} U_s \\ 0\end{pmatrix}\) of the (MHD) system is nonlinearly unstable in the following sense:
	
  There exists \(T = (T_1,T_2,T_3)\in(\mb N^*)^3\) and initial values 
  \begin{equation}
    \begin{pmatrix} u_0 \\ b_0\end{pmatrix}\in\mc H_T^s,
  \end{equation}
    and \(C_0>0\) such that for all \(\delta>0\), the solution \(\begin{pmatrix} u_\delta \\ b_\delta \end{pmatrix}\) of \eqref{MHDPert} with initial value \(\delta\begin{pmatrix} u_0 \\ b_0\end{pmatrix}\) satisfies
  \begin{equation}
    \left\|\begin{pmatrix} u_\delta \\ b_\delta \end{pmatrix}(t_\delta)\right\|_{\mc H_T^s}\ge C_0
  \end{equation}
  for some time \(t_\delta\). Furthermore, we have the estimate on \(t_\delta\) as \(\delta\rightarrow0\):
  \begin{equation}
    t_\delta\sim_{\delta\rightarrow0} -\frac{\ln\delta}{C}
  \end{equation}
\end{theoreme}
As expected, this theorem expresses the instability of the MHD system, up to consider a larger periodic box \([0,T_1]\times[0,T_2]\times[0,T_3]\). Let us point out that \((T_1,T_2,T_3)\in(\mb N^*)^3\), so that \(U_s\) is also periodic on this larger scale.

In particular, for any
\begin{equation}
  \begin{pmatrix}
    u_0\\
    b_0
    \end{pmatrix}\in\mc H_T^s,
  \end{equation}
and any \(\delta>0\), system \eqref{MHDPert} has a unique solution
\begin{equation} 
  \begin{pmatrix}
    u_\delta\\
    b_\delta
    \end{pmatrix}\in\mc C^0([0,T^*],H_T^s)
  \end{equation}
for some \(T^*>0\) with initial data \(\delta\begin{pmatrix} u_0,b_0\end{pmatrix}\). Moreover, for \(\delta>0\) small enough, it is well-known that one can take \(T^*=+\infty\).

Up to our knowledge, our theorem is the first full justification of the \(\alpha\)-effect, in the nonlinear context.

In particular, it extends substentially the article \cite{gerardvaret2006osi} by David G\'erard-Varet. Recast in our variables, this article considers the case of large magnetic diffusion, that is \(R_m=\varepsilon\ll1\). Briefly, it is shown that for all \(m\in\mb N\), there are solutions of \eqref{MHDPert} (on a larger box depending on \(\varepsilon\)) that go from amplitude \(\varepsilon^m\) initially to an amplitude \(\eta = \eta(m)>0\) independent of \(\varepsilon\) (but depending on \(m\)).

Note that, as \(\eta\) depends on \(m\), this result does not yield Liapunov instability at fixed small \(\varepsilon\). Moreover, the study in \cite{gerardvaret2006osi} corresponds to vanishing magnetic Reynolds number. In our theorem, the threshold \(R_m^0\) is of order \(1\), which is much more satisfactory from a physical point of view.

The outline of the paper is as follow:

The first and main part of the paper is the constructino of an exact unstable eigenmode for the linearization of \eqref{MHDPert}. It involves sharp spectral analysis arguments. Then, on the basis of this linear analysis, the proof of nonlinear instability is given in section 3.


\section{Linear instability}
In this section we will establish linear instability for the induction part of the MHD system. In other words, we will exhibit an unstable eigenmode of the linear operator
\begin{equation}
	\fonction[][b][\nabla\wedge(U\wedge b) - \frac1{R_m}\Delta b]
	\end{equation}
for a large class of \(U\). This will be possible, through an accurate mathematical analysis of the \emph{alpha effect} described in the introduction.

Let us write down the induction equation:
\begin{equation}
	\label{MHD_induction_ne}
	\begin{split}
	\partial_tb - \nabla\wedge(U\wedge b) - \frac1{R_m}\Delta b = 0,\\
	\nabla\cdot b = 0,
	\end{split}
	\end{equation}
where \(U\in\mc H_0^\infty(\mb T^3)\) is independent of time.

In this section we prove the following linear instability result:
\begin{theoreme}
	\label{theoreme_lin}
 	Let \(s\ge0\) be a real number. There exists a dense subset \(\ms P\) of the unit ball of \(\mc H_0^\infty(\mb T^3)\) and a positive real \(R_m^0\) such that for all positive \(R_m<R_m^0\), and all \(U\in\ms P\), there exists \(\xi\in\mb (2\pi\mb Q^*)^3\) and a number \(\varepsilon_0>0\) such that for all \(0<\varepsilon<\varepsilon_0\), there exists a solution \(B\) of \eqref{MHD_induction_ne} of the form
\begin{equation}
  b(t,\theta) = e^{\iim\varepsilon\xi\cdot\theta}e^{\lambda^\varepsilon t}b^\varepsilon(\theta)
  \end{equation}
where \(b^\varepsilon\in\mc H^s(\mb T^3)\), with exponential growth when \(t\rightarrow+\infty\) (that is \(\Re(\lambda^\varepsilon)>0\)).
 	\end{theoreme}

\subsection{Linearized equation}
Let \(U\in\mc H_0^\infty(\mb T^3)\) be a spatially periodic velocity field in \(\mb R^3\) (of period \(1\)).

Since we will look for the growth of the magnetic field on a larger space-scale \(1/\varepsilon\), we will directly look for a solution \(b\) under the form
\begin{equation}
	b(\theta,t) = B^\varepsilon(\varepsilon \theta,\theta,\varepsilon t) = \bar B^\varepsilon(\varepsilon\theta,\varepsilon t) + \tilde B^\varepsilon(\varepsilon\theta,\theta,\varepsilon t),
	\end{equation}
where \(B^\varepsilon = B^\varepsilon(x,\theta,t)\) depends both on the original variable \(\theta\) and a slow variable \(x = \varepsilon\theta\). One can decompose \(B^\varepsilon\) as the sum of a large-scale field \(\bar B^\varepsilon\) which doesn't depend on the smaller variable, and a small-scale zero-mean-value field \(\tilde B^\varepsilon\), i.e. \(\int_{\mb T^3}\tilde B^\varepsilon(x,\theta,t)\ud\theta = 0\), as suggested in the introduction.

In view of \eqref{MHD}, the new field \(B^\varepsilon\) should satisfy (dropping the prime over time variable \(t\)):
\begin{equation}
	\label{MHD_induction}
	\left\{\begin{array}{l}
		\partial_t B^\varepsilon - \nabla_x\wedge(U\wedge B^\varepsilon) - \frac1\varepsilon\nabla_\theta\wedge(U\wedge B^\varepsilon) - \frac{\varepsilon}{R_m}\Delta_xB^\varepsilon - \frac{1}{\varepsilon R_m}\Delta_\theta B^\varepsilon - \frac{2}{R_m}\sum_{i=1}^d\partial_{x_i\theta_i}B^\varepsilon = 0,\\
		\nabla_x\cdot B^\varepsilon + \frac1\varepsilon\nabla_\theta\cdot B^\varepsilon = 0.
		\end{array}\right.
	\end{equation}

The idea here is to find an exact solution of this new equation in \((t,x,\theta)\) (or more precisely to prove that there exists one that fits our needs). We will find this solution as a Fourier mode in \(x\):
\begin{equation}
	B^\varepsilon = e^{\iim\xi\cdot x}{B^\varepsilon}'(t,\theta).
	\end{equation}
Dropping the prime again, we get:
\begin{equation}
	\label{MHD_induction_Fourier}
	\left\{\begin{array}{l}
		\partial_t B^\varepsilon - \iim\xi\wedge(U\wedge B^\varepsilon) - \frac1\varepsilon\nabla_\theta\wedge(U\wedge B^\varepsilon) + \frac{\varepsilon}{R_m}\xi^2B^\varepsilon - \frac{1}{\varepsilon R_m}\Delta_\theta B^\varepsilon - \frac{2}{R_m}\iim\sum_{i=1}^d\xi_i\partial_{\theta_i}B^\varepsilon = 0,\\
		\iim\xi\cdot B^\varepsilon + \frac1\varepsilon\nabla_\theta\cdot B^\varepsilon = 0.
		\end{array}\right.
	\end{equation}
We denote by \(A^\varepsilon = A^0 + \varepsilon A^1 + \varepsilon^2A^2\) the operator associated to the first equation of \eqref{MHD_induction_Fourier}, that is
\begin{equation}
 	\begin{split}
	A^0 & = \nabla_\theta\wedge(U\wedge \cdot) + \frac{1}{R_m}\Delta_\theta \\
        A^1 & = \iim\xi\wedge(U\wedge \cdot) + \frac{2}{R_m}\iim\sum_{i=1}^d\xi_i\partial_{\theta_i} \\
	A^2 & = -\frac{1}{R_m}\xi^2\Id,\\
        \end{split}
	\end{equation}
where \(\fonction{A^\varepsilon}{\mc H^{s+2}(\mb T^3)}{\mc H^s(\mb T^3)}\),
such that the first equation in \eqref{MHD_induction_Fourier} writes 
\begin{equation}
	\partial_t B^\varepsilon = \frac1\varepsilon A^\varepsilon B^\varepsilon.
	\end{equation}
We want to find a solution in power series of \(\varepsilon\) of that equation of type
\begin{equation}
	B^\varepsilon(\theta,t) = e^{\lambda^\varepsilon t}B_0^\varepsilon(\theta).
	\end{equation}
Thus, the equation writes:
\begin{equation}
	\label{equation_operateur}
	A^\varepsilon B_0^\varepsilon = \varepsilon\lambda^\varepsilon B_0^\varepsilon,
	\end{equation}
which is an eigenvalue problem. We will split this equation in powers of \(\varepsilon\) to obtain \(3\) equations.

\subsection{Study of the first term in the development}
The first equation, corresponding to \(\varepsilon = 0\), writes
\begin{equation}
	A^0B_0^0 = 0,
	\end{equation}
that is \(\frac1{R_m}\Delta_\theta B_0^0 + \nabla_\theta\wedge(U\wedge B_0^0) = 0\), which depends only on the independent variable \(\theta\).

We will look for solution of type \(B_0^0(\theta) = \bar B_0^0 + \tilde B_0^0(\theta)\) as proposed before:
\begin{equation}
	\label{equation_L_rm}
	\frac1{R_m}\Delta_\theta \tilde B_0^0 + \nabla_\theta\wedge(U\wedge \tilde B_0^0) = - \nabla_\theta\wedge(U\wedge \bar B_0^0).
	\end{equation}
Since the Laplacian
\begin{equation}
	\fonction{\Delta_\theta}{\mc H_0^{s+2}(\mb T^3)\subset\mc H_0^s(\mb T^3)}{\mc H_0^s(\mb T^3)}
	\end{equation}
is an unbounded invertible operator, one can introduce
\begin{equation}
	\fonction{\Delta_\theta^{-1}\nabla_\theta\wedge(U\wedge\cdot)}{\mc H_0^s(\mb T^3)}{\mc H_0^{s+1}(\mb T^3)\subset \mc H_0^s(\mb T^3)}
	\end{equation}
which is a compact operator of domain \(\mc H_0^s(\mb T^3)\).

Thus, by classical result on the spectrum of compact operators, the operator
\begin{equation}
	L_{R_m} = \frac1{R_m}\Delta_\theta \cdot + \nabla_\theta\wedge(U\wedge \cdot) = \Delta_\theta(\frac1{R_m}\Id + \Delta_\theta^{-1}\nabla_\theta\wedge(U\wedge \cdot))
	\end{equation}
is invertible as an operator from \(\mc H_0^{s+2}(\mb T^3)\) to \(\mc H_0^s(\mb T^3)\) but for some countable set \(C_L\) of \(R_m\) which can be written as a sequence of numbers growing to \(+\infty\). In particular it is invertible for \(R_m\) smaller than a given value \(R_m^0\).

Since \(-\nabla_\theta\wedge(U\wedge\bar B_0^0)\in\mc H_0^s(\mb T^3)\) for any \(\bar B_0^0\in\mb R^3\), the equation \eqref{equation_L_rm} admits exactly one solution for every \(\bar B_0^0\in\mb R^3\). Then for any \(R_m\) out of \(C_L\), the equation
\begin{equation}
	L_{R_m}\tilde B_0^0 = - \nabla_\theta\wedge(U\wedge \bar B_0^0)
	\end{equation}
is equivalent to
\begin{equation}
	\tilde B_0^0 = \ms L(\theta)\bar B_0^0.
	\end{equation}

where \(\ms L(\theta)\) is the operator from \(\mb R^3\) to \(\mc H_0^{s+2}(\mb T^3)\) defined by
\begin{equation}
	\ms L(\theta) = - L_{R_m}^{-1}\nabla_\theta\wedge(U\wedge \cdot)
	\end{equation}

That means that the \(A^0\) operator has a kernel of dimension \(3\), \(\fonction{A^0}{\mc H^{s+2}(\mb T^3)\subset \mc H^s(\mb T^3)}{\mc H^s(\mb T^3)}\), that is \(\ker A^0 = \{\bar B+\ms L(\theta)\bar B,\bar B\in\mb R^3\}\).

Using Lax-Milgram theorem and elliptic regularity, we know that the operator
\begin{equation}
	\frac1{R_m}\Delta + \nabla_\theta\wedge(U\wedge\cdot) - \lambda\Id = A^0-\lambda\Id
	\end{equation}
is invertible for large enough \(\lambda\) if seen as an operator from \(\mc H_0^{s+2}\) to \(\mc H_0^s\). Thus the resolvant of \(A^0\) is compact in at least one \(\lambda\), and all spectral projections are of finite dimension, that is \(A^0\) has eigenvalues of finite dimension, isolated and they are the only spectral values of \(A^0\) (see \cite{kato1995perturbation} p.187).

Let \(\sigma\) be a smooth path surrounding \(0\) and no other eigenvalue of \(A^0\). For \(\varepsilon\) sufficiently small, since \(A^\varepsilon\) varies continuously on \(\varepsilon\), \(\sigma\) doesn't cross the spectrum of \(A^\varepsilon\), and the operator
\begin{equation}
	P^\varepsilon := \frac1{2\iim\pi}\int_\sigma(A^\varepsilon - \lambda\Id)^{-1}\ud\lambda.
	\end{equation}
is the spectral projector associated to the eigenvalues of \(A^\varepsilon\) which tend to \(0\) as \(\varepsilon\) tends to \(0\). (\(0\) is an eigenvalue of degree \(3\) in the case of \(A^0\), but this eigenvalue will vary -- and maybe split -- around \(0\) for positive \(\varepsilon\), without crossing \(\sigma\) if \(\varepsilon\) is small enough). Note that the path depends on the choice of \(\xi\).

Morover, still for \(\varepsilon\) small enough, \(P^\varepsilon\) also varies analytically on \(\varepsilon\).

Furthermore, by definition, \(P^0\) is the spectral projection of \(A^0\) on its eigenvalue \(0\), and commutes with \(A^0\) (\cite{kato1995perturbation} again, thm 6.17 p 178).

On \(\im P^0\), \(A^0\) has then only \(0\) as an eigenvalue. Since \(\im P^0\) is of finite dimension, \(A^0\) is nilpotent on this same space.
Until the end of the subsection, we define \(\mc A^0 = {A^0}_{|\im P^0}\).

Let \(n\) be the integer such that
\begin{equation}
	\ker \mc A^0\subset \ker (\mc A^0)^2 \subset \cdots \subset \ker (\mc A^0)^n = \im P^0,
	\end{equation}

We will now prove that \(n = 1\), that is \(\ker \mc A^0 = \ker (\mc A^0)^2\):\\
Let \(B\in\ker (\mc A^0)^2\). Then we have the equivalences:
\begin{equation}
	\begin{split}
		(\mc A^0)^2B = 0 &\iff \mc A^0(\widetilde{\mc A^0B}) + \mc A^0(\underbrace{\overline{\mc A^0B}}_{=0}) = 0, \\
			&\iff \widetilde{\mc A^0B} = \bar{B'}+\ms L(\theta)\bar{B'}\text{ and }\overline{\mc A^0B} = 0, \\
			&\iff \widetilde{\mc A^0B} = 0 \text{ and }\overline{\mc A^0B} = 0, \\
			&\iff \mc A^0B = 0.
		\end{split}
	\end{equation}
Thus \(\ker \mc A^0 = \ker (\mc A^0)^2 = \im P^0\).

\subsection{Simplification of the equation}
Back to \eqref{equation_operateur}, we look for solutions of that equation in \(\im P^\varepsilon\), that is in the sum of the proper-spaces corresponding to the spectral values of \(A^\varepsilon\) surrounded by \(\sigma\):
\begin{equation}
	P^\varepsilon A^\varepsilon P^\varepsilon B^\varepsilon = \varepsilon\lambda^\varepsilon B^\varepsilon.
	\end{equation}
Developping \(P^\varepsilon\) as the beginning of a series in \(\varepsilon\), \(P^\varepsilon = P^0 + \varepsilon P^1 + \varepsilon^2 {P'}^\varepsilon\), we have
\begin{equation}
	P^\varepsilon A^\varepsilon P^\varepsilon = P^0A^0P^0 + \varepsilon P^1A^0P^0 + \varepsilon P^0A^0P^1 + \varepsilon P^0A^1P^0 + \varepsilon^2{A'}^\varepsilon,
	\end{equation}
where both \({P'}^\varepsilon\) and \({A'}^\varepsilon\) depend analytically on \(\varepsilon\).

Since \(P^0\) and \(A^0\) commutes and \(A^0P^0 = 0\) (see results of previous subsection), the first three terms vanish, thus leaving us with the equation:
\begin{equation}
	\label{equation_projetee}
	P^0 A^1 P^0 B^\varepsilon +\varepsilon {A'}^\varepsilon B^\varepsilon = \lambda^\varepsilon B^\varepsilon.
	\end{equation}

In the next subsections, we look for a solution of this equation for \(\varepsilon = 0\). We shall then conclude using the continuity of eigenvalues as functions of \(\varepsilon\).

\subsection{Case \texorpdfstring{$\varepsilon=0$}{epsilon = 0}}
The first order in \(\varepsilon\) of the equation \eqref{equation_projetee} is equivalent to the system
\begin{equation}
	\left\{\begin{array}{l}
		A^1B = \lambda B \\
		\tilde B = \ms L(\theta)\bar B \\
		\end{array}\right.
	\end{equation}
where \(A^1 = \iim\xi\wedge(U\wedge \cdot) - \frac{2\iim}{R_m}\sum_{i=1}^d\xi_i\partial_{\theta_i}\).
The mean part of the first equation writes
\begin{equation}
	\iim\xi\wedge(\overline{U\wedge B}) = \lambda\bar B,
	\end{equation}
that is
\begin{equation}
	\iim\xi\wedge(\overline{U\wedge\ms L(\theta)\bar B}) + \underbrace{\iim\xi\wedge(\overline{U\wedge\bar B})}_{=0} = \lambda\bar B,
	\end{equation}
\begin{equation}
	\label{alpha_eq}
	\iim\xi\wedge(\alpha\bar B) = \lambda\bar B,
	\end{equation}
where \(\alpha\) is a constant matrix (that depends on \(U\) and \(R_m\)),
\begin{equation}
	\alpha b = \int_{\theta\in\mb T^3}U(\theta)\wedge\ms L(\theta)b,
	\end{equation}
which we will study in the next subsection.

Thus we now have to find the eigenvalues of the operator \(\iim\xi\wedge(\alpha \cdot)\), that is to find the eigenvalues of the matrix \(A^\xi\alpha\), where
\begin{equation}
	A^\xi = \begin{pmatrix}
	0 & -\iim\xi_3 & \iim\xi_2 \\
	\iim\xi_3 & 0 & -\iim\xi_1 \\
	-\iim\xi_2 & \iim\xi_1 & 0 \\
	\end{pmatrix}.
	\end{equation}

\subsection{Study of the matrix \texorpdfstring{$\alpha$}{alpha}}
The \(\alpha\) matrix is defined for any \(b\in\mb R^3\) by
\begin{equation}
	\alpha b = \int_{\theta\in\mb T^3}U(\theta)\wedge\ms L(\theta)b = -\int_{\theta\in\mb T^3}U(\theta)\wedge\left(L_{R_m}^{-1}\nabla_\theta\wedge(U(\theta)\wedge b)\right).
	\end{equation}

We shall write this matrix as the sum of his antisymmetric part \(\alpha^A\) and symmetric part \(\alpha^S\) (which are both real).

Now we compute the eigenvalues of those matrices:\\
Since \(\alpha^A\) is antisymmetric, his eigenvalues are \(0,\pm\iim\|V_U\|\) where \(V_U\) is a vector of \(\mb R^3\) such that \(\forall b\in\mb R^3,\alpha^Ab = V_U\wedge b\).\\

For \(\alpha^S\), we will apply a slight variation to \(U\) so that the matrix changes to one of distinct non-zero eigenvalues (this doesn't change what we said about the antisymmetric part).

\subsubsection{Symmetry of \texorpdfstring{$\alpha$}{alpha}}
We have in the expression of \(\alpha\):
\begin{equation}
	-L_{R_m}^{-1}\nabla_\theta\wedge(U\wedge b) = \left(\frac1{R_m}\Id - \ms A\right)^{-1}\big(-\Delta_\theta^{-1}\nabla_\theta\wedge(U\wedge b)\big) = \left(\frac1{R_m}\Id - \ms A\right)^{-1}\ms Ab,
	\end{equation}
where
\begin{equation}
	\ms A = \fonction{-\Delta_\theta^{-1}\nabla_\theta\wedge(U\wedge \cdot)}{\mc H_0^{s+2}}{\mc H_0^{s+3}\subset\mc H_0^{s+2}}. \label{A.def}
	\end{equation}
For small enough \(R_m\) (such that \(1/R_m\) is bigger than the highest eigenvalue of the compact operator \(\ms A\), say \(R_m<2R_m^0\) where \(R_m^0>0\) -- we will justify this factor \(2\) as a precaution at the end of the subsection), we can express the resolvent as a series:
\begin{equation}
	-L_{R_m}^{-1}\nabla_\theta\wedge(U\wedge b) = \sum_{n=0}^{+\infty}R_m^{n+1}\ms A^{n+1}b = \sum_{n=1}^{+\infty}R_m^n\ms A^nb.
	\end{equation}
Furthermore, \eqref{A.def} gives in Fourier series, for \(k\neq0\) and \(f\in\mc H^s(\mb T^3)\):
\begin{equation}
	\widehat{\ms Af}(k) = -\frac1{2\pi}|k|^{-2}(-\iim k)\wedge\left(\sum_{k'}\hat U(k')\wedge\hat f(k-k')\right).
	\end{equation}
Note that, by definition of \(\Delta_\theta^{-1}\), \(\widehat{Af}(0) = 0\).

thus
\begin{equation}
	\begin{split}
		\alpha b &= \int_{\theta\in\mb T^3}U(\theta)\wedge\sum_{n=1}^{+\infty}(R_m\ms A)^nb \\
		& = \sum_{n=1}^{+\infty}R_m^n\underbrace{\int_{\theta\in\mb T^3}U(\theta)\wedge\ms A^nb}_{\widehat{U\wedge\ms A^nb}(0)}. \\
	\end{split}
	\end{equation}
Recursively, we can compute \(\widehat{U\wedge\ms A^nb}(0)\) and
\begin{equation}
	\alpha b = \sum_{n=1}^{\infty}R_m^n\sum_{\sum_{j=1}^{n+1}k_j = 0}\hat U(k_{n+1})\wedge\left(\frac{\iim m_{n}}{2\pi|m_{n}|^2}\wedge\left(\hat U(k_{n})\wedge\cdots\left(\frac{\iim m_1}{2\pi|m_1|^2}\wedge\left(\hat U(k_1)\wedge b\right)\right)\cdots\right)\right),
	\end{equation}
where the second sum is actually restricted to tuples where all partial sums \(m_i=\sum_{j=1}^ik_j\) are non-zero (otherwise there is a term of the form \(\widehat{\ms A^kb}(0)\) which is zero by definition of \(\ms A\)).

This last equality writes:
\begin{equation}
	\alpha b = \sum_{n=1}^{\infty}R_m^n\alpha^{(n+1)}b,
	\end{equation}
where
\begin{equation}
	\alpha^{(n)} = \sum_{\sum_{j=1}^{n}k_j = 0}\hat U(k_n)\wedge\left(\frac{\iim m_{n-1}}{2\pi|m_{n-1}|^2}\wedge\left(\hat U(k_{n-1})\wedge\cdots\left(\frac{\iim m_1}{2\pi|m_1|^2}\wedge\left(\hat U(k_1)\wedge b\right)\right)\cdots\right)\right),
	\end{equation}
and \(\alpha^{(n)}\) have the same symmetry as \(n\) (by recurrence, or see \cite{childress1967construction}).

\subsubsection{Eigenvalues of \texorpdfstring{$\alpha$}{alpha}}
The first symmetric matrix in power series is \(\alpha^{(2)}\)~:\\
Let \(U\in\mc H_0^\infty(\mb T^3)\), and \(\eta > 0\). We want \(\tilde U\) such that \(\alpha^{(2)}(\tilde U)\) is of distinct non-zero eigenvalues and \(d(U,\tilde U)<\eta\), where \(d\) is the usual distance induced by semi-norms on \(\ms H_0^\infty(\mb T^3)\): we define now
\begin{equation}
	U^j(\theta) = \sum_{|k|\le j}\hat U(k)e^{2\iim k\cdot\theta\pi}.
	\end{equation}
For large enough \(j\), we have \(d(U^j,U)\le \eta/2\). Denote then \(\tilde U = U^j + \sum_{i=1}^3\delta_iV^i\) where \(\delta_i>0\) is determined later, and
\begin{align}
	V^1 &= \begin{pmatrix}
		\sin(2\pi(j+1)\theta_3) + \cos(2\pi(j+1)\theta_2)\\ \cos(2\pi(j+1)\theta_3) \\ \sin(2\pi(j+1)\theta_2)
		\end{pmatrix},&
	V^2 &= \begin{pmatrix}
		\sin(2\pi(j+2)\theta_3) \\ \sin(2\pi(j+2)\theta_1) + \cos(2\pi(j+2)\theta_3)\\ \cos(2\pi(j+2)\theta_1)
		\end{pmatrix}, \nonumber\\
	V^3 &= \begin{pmatrix}
		\cos(2\pi(j+3)\theta_2) \\ \sin(2\pi(j+3)\theta_1) \\ \sin(2\pi(j+3)\theta_2) + \cos(2\pi(j+3)\theta_1)
		\end{pmatrix}.
	\end{align}
Then:
\begin{equation}
	\begin{split}
	\hat V^1(k)_y &= \int e^{-2\iim k\cdot\theta\pi}V_y^1(\theta)\ud\theta = \left\{\begin{array}{l} 0\text{ if }k_1\neq 0\text{ or }k_2\neq 0 \\
	\int_0^1e^{-2\iim\pi k_3\theta_3}\cos(2\pi(j+1)\theta_3)\ud\theta_3\text{ otherwise}\end{array}\right. \\
	& = \left\{\begin{array}{l} 0\text{ if }(k_3\neq j+1\text{ and }k_3\neq -j-1)\text{ or }k_1 \neq 0\text{ or }k_2 \neq 0 \\
	\frac12\text{ otherwise}\end{array}\right.
	\end{split}.
	\end{equation}
Thus
\begin{equation}
	\hat V^1(k) = \begin{pmatrix}
		\frac12\delta_{k_1=0}\delta_{k_3=0}(\delta_{k_2 = j+1} + \delta_{k_2 = -j-1}) +
			\frac\iim2\delta_{k_1=0}\delta_{k_2=0}(\delta_{k_3 =-j-1} - \delta_{k_3 = j+1}) \\
		\frac12\delta_{k_1=0}\delta_{k_2=0}(\delta_{k_3 = j+1} + \delta_{k_3 = -j-1}) \\
		\frac\iim2\delta_{k_1=0}\delta_{k_3=0}(\delta_{k_2 =-j-1} - \delta_{k_2 = j+1}) \\
		\end{pmatrix},
	\end{equation}
and similarily for the other ones:
\begin{equation}
	\hat V^2(k) = \begin{pmatrix}
		\frac\iim2\delta_{k_1=0}\delta_{k_2=0}(\delta_{k_3 =-j-2} - \delta_{k_3 = j+2}) \\
		\frac12\delta_{k_1=0}\delta_{k_2=0}(\delta_{k_3 = j+2} + \delta_{k_3 = -j-2}) +
			\frac\iim2\delta_{k_2=0}\delta_{k_3=0}(\delta_{k_1 =-j-2} - \delta_{k_1 = j+2}) \\
		\frac12\delta_{k_2=0}\delta_{k_3=0}(\delta_{k_1 = j+2} + \delta_{k_1 = -j-2}) \\
		\end{pmatrix},
	\end{equation}
\begin{equation}
	\hat V^3(k) = \begin{pmatrix}
		\frac12\delta_{k_1=0}\delta_{k_3=0}(\delta_{k_2 = j+3} + \delta_{k_2 = -j-3}) \\
		\frac\iim2\delta_{k_2=0}\delta_{k_3=0}(\delta_{k_1 =-j-3} - \delta_{k_1 = j+3}) \\
		\frac12\delta_{k_2=0}\delta_{k_3=0}(\delta_{k_1 = j+3} + \delta_{k_1 = -j-3}) +
			\frac\iim2\delta_{k_1=0}\delta_{k_3=0}(\delta_{k_2 =-j-3} - \delta_{k_2 = j+3}) \\
		\end{pmatrix}.
	\end{equation}
For \(V^i(k)\) to be non-zero, we need that \(k_i = 0\) and \(k_l\in\{0,j+i,-j-i\}\) when \(l\neq i\).

Then \(\hat V^1(-k)\wedge\left(\iim\frac k{|k|^2}\wedge\left(\hat V^2(k)\wedge B\right)\right) = 0\) as soon as \(k_1\neq0\) or \(k_2\neq0\). If both are null, then \(\hat V^1(k) = 0\) but if \(k_3\in\{0,j+1,-j-1\}\), and \(\hat V^2(k) = 0\) but if \(k_3\in\{0,j+2,-j-2\}\). Thus the product is always zero.

Only remains now non-crossed terms:
\begin{equation}
	\hat V^1(-k)\wedge\left(\iim\frac k{|k|^2}\wedge\left(\hat V^1(k)\wedge B\right)\right) = \hat V^1(-k)\wedge\hat V^1(k).\left(\iim\frac k{|k|^2}\cdot B\right).
	\end{equation}
For it to be non-zero, we need that \(k_1 = 0\). If furthermore \(k_2\neq0\) and \(k_3\neq0\), then \(\hat V^1(k) = 0\).
\begin{center}
\begin{tabular}{c|c}
	If \(k_3 = 0\) and \(k_2\neq0\) & If \(k_3\neq0\) and \(k_2=0\) \\
	\(\hat V^1(k) = \begin{pmatrix} \frac12(\delta_{k_2=j+1} + \delta_{k_2=-j-1}) \\ 0 \\ \frac\iim2(\delta_{k_2=-j-1} - \delta_{k_2=j+1}) \end{pmatrix}\) & 
	\(\hat V^1(k) = \begin{pmatrix} \frac\iim2(\delta_{k_3=-j-1} - \delta_{k_3=j+1}) \\ \frac12(\delta_{k_3=j+1} + \delta_{k_3=-j-1}) \\ 0 \end{pmatrix}\)\\
	Thus for \(k_2 = j+1\),
	\(\hat V^1(-k)\wedge\hat V^1(k) = \frac12\begin{pmatrix} 0 \\ \iim \\ 0\end{pmatrix}\) & 
	Thus for \(k_3 = j+1\),
	\(\hat V^1(-k)\wedge\hat V^1(k) = \frac12\begin{pmatrix} 0 \\ 0 \\ \iim\end{pmatrix}\)
	\end{tabular}
\end{center}
and thus
\begin{equation}
	\begin{split}
	\sum_{k\neq0} &\hat V^1(-k)\wedge\hat V^1(k)\left(\iim\frac k{|k|^2}\cdot B\right) \\
	& = \sum_{k_3=k_1=0}\hat V^1(-k)\wedge\hat V^1(k)\left(\iim\frac k{|k|^2}\cdot B\right) + \sum_{k_2=k_1=0}\hat V^1(-k)\wedge\hat V^1(k)\left(\iim\frac k{|k|^2}\cdot B\right) \\
	& = \frac12\begin{pmatrix} 0 \\ \iim \\ 0\end{pmatrix}.\left(\frac \iim{|j+1|^2}\begin{pmatrix} 0 \\ j+1 \\ 0\end{pmatrix}\cdot B\right) + \frac12\begin{pmatrix} 0 \\ -\iim \\ 0\end{pmatrix}.\left(\frac \iim{|j+1|^2}\begin{pmatrix} 0 \\ -j-1 \\ 0\end{pmatrix}\cdot B\right) \\
	&\quad + \frac12\begin{pmatrix} 0 \\ 0 \\ \iim \end{pmatrix}.\left(\frac \iim{|j+1|^2}\begin{pmatrix} 0 \\ 0 \\j+1 \end{pmatrix}\cdot B\right) + \frac12\begin{pmatrix} 0 \\ 0 \\ -\iim \end{pmatrix}.\left(\frac \iim{|j+1|^2}\begin{pmatrix} 0 \\ 0 \\ -j-1\end{pmatrix}\cdot B\right) \\
	& = -\frac1{j+1}\begin{pmatrix} 0 \\ B_y \\ B_z\end{pmatrix}.
	\end{split}
	\end{equation}
We do the same for the other ones:
\begin{equation}
	\sum_{k\neq0} \hat V^2(-k)\wedge\hat V^2(k)\left(\iim\frac k{|k|^2}\cdot B\right) = \frac{-1}{j+2}\begin{pmatrix} B_x \\ 0 \\ B_z\end{pmatrix},\quad\sum_{k\neq0} \hat V^3(-k)\wedge\hat V^3(k)\left(\iim\frac k{|k|^2}\cdot B\right) = \frac{-1}{j+3}\begin{pmatrix} B_x \\ B_y \\ 0\end{pmatrix}.
	\end{equation}
Finally, for all \(i\in\lb1,3\rb\), since the spectrum of \(U^j\) and \(V^i\) are disjoint:
\begin{equation}
	\sum_{k\neq0} \hat U^j(-k)\wedge\hat V^i(k)\left(\iim\frac k{|k|^2}\cdot B\right) = 0,
	\end{equation}
and then
\begin{equation}
	\alpha^{(2)}(\tilde U) = \alpha^{(2)}(U^j) - \frac{\delta_1}{j+1}\begin{pmatrix} 0 & 0 & 0 \\ 0 & 1 & 0 \\ 0 & 0 & 1\end{pmatrix} - \frac{\delta_2}{j+2}\begin{pmatrix} 1 & 0 & 0 \\ 0 & 0 & 0 \\ 0 & 0 & 1\end{pmatrix} - \frac{\delta_3}{j+3}\begin{pmatrix} 1 & 0 & 0 \\ 0 & 1 & 0 \\ 0 & 0 & 0\end{pmatrix},
	\end{equation}
where the \(\delta_i\) can be chosen as small as we want. Thus \(\alpha^{(2)}(U)\) is as near as we want from a matrix of distinct non-zero eigenvalues, and it's the same for \(\alpha^S(U)\) for \(R_m<R_m^0\) (NB: the set in which \(R_m\) can be taken will vary continuously with \(U\), thus the precaution we took at the beginning of the subsection to be sufficiently far from the limit)


\subsubsection{End of the proof}
In the sequel of this subsection, we now work with \(\tilde U\) instead of \(U\), droping the tilda. Thus, \(\alpha^S\) has distinct non-zero eigenvalues. In Fourier transform, we now have to find the eigenvalues of the matrix \(A^\xi\alpha^A + A^\xi\alpha^S\), where
\begin{equation}
	A^\xi = \iim\begin{pmatrix} 0 & -\xi_3 & \xi_2 \\ \xi_3 & 0 & -\xi_1 \\ -\xi_2 & \xi_1 & 0 \end{pmatrix}.
	\end{equation}
Let there be \(P\in\ms O_3(\mb R)\) such that
\begin{equation}
	{}^tP\alpha^SP = \begin{pmatrix} 
		\alpha_1 & 0 & 0 \\
		0 & \alpha_2 & 0\\
		0 & 0 & \alpha_3\end{pmatrix}
		= D,
	\end{equation}
and write \(\zeta = {}^tP\xi\). The matrix \(A^\zeta D\) has \(\lambda_0 = 0\) eigenvalue, and two opposite eigenvalues given by
\begin{equation}
	\lambda_\pm = \pm\sqrt{\zeta_1^2\alpha_2\alpha_3 + \zeta_2^2\alpha_3\alpha_1 + \zeta_3^2\alpha_1\alpha_2}.
	\end{equation}
Those eigenvalues are real for any \(\zeta\) if all eigenvalues of \(\alpha^S\) have the same sign. Otherwise, there exists a cone with vertex \(0\) outside of which they are. Furthermore, for \(\zeta\) in the complement \(\Omega\) of this cone, \(\lambda_+\) stays positive, and \(\lambda_-\) stays negative. Thus, the eigenvectors of \(A^\zeta D\) depend regularly of \(\zeta\). From now on, we fix some \(\zeta\) in this good set \(\Omega\). We can also assume that \(\xi = P\zeta\in(2\pi\mb Q_+^*)^3\).

Let \(\beta\) be an eigenvector associated to \(\lambda_+\). Then \(A^\zeta D\beta = \iim\zeta\wedge (D\beta) = \lambda_+\beta\), so \(\zeta\cdot\beta = 0\). Thus, we have \(A^\zeta{}^tP\alpha^AP\beta = \zeta\wedge(\gamma\wedge\beta) = -\iim(\zeta\cdot\gamma)\beta\), where \(\gamma\) is a vector such that \(\forall b\in\mb R^3,{}^tP\alpha^APb = \gamma\wedge b\). Thus \(\beta\) is also an eigenvector of \(A^\zeta{}^tP\alpha^AP\), associated with the eigenvalue \(-\iim(\zeta\cdot\gamma)\).

Now changing back to \(\xi = P\zeta\), we have that
\begin{equation}
	A^\zeta = {}^tPA^\xi P,
	\end{equation}
so that \(\beta\) is an eigenvector of \({}^tPA^\xi\alpha^AP+{}^tPA^\xi PD\) associated to the eigenvalue \(\lambda_+ -\iim(\zeta\cdot\gamma) = \lambda_+ -\iim(\xi\cdot P\gamma)\). And thus \(P\beta\) is an eigenvector of \(A^\xi\alpha^A + A^\xi\alpha^S\) associated to an eigenvalue of positive real part, that is a growing mode solution of the first term of \eqref{equation_projetee}.

\subsection{Conclusion to the linear case}
We have proven now the existence of an eigenvector \(b^0_0\) of the equation
\begin{equation}
	P^0 A^1 P^0 B^\varepsilon +\varepsilon {A'}^\varepsilon B^\varepsilon = \lambda^\varepsilon B^\varepsilon
	\end{equation}
for \(\varepsilon = 0\) with positive real part, that is an exponentially-growing mode of the form
\begin{equation}
	B^\varepsilon(x,\theta,t) = b^0_0(\theta,\zeta) e^{(\lambda_+-\iim\zeta\cdot\gamma)t + \iim\zeta\cdot x}
	\end{equation}
of the induction equation for given velocity field \(U_s\) and small enough \(R_m\) (and \(\varepsilon = 0\)). Since the operator \(\frac1\varepsilon P^\varepsilon A^\varepsilon P^\varepsilon\) varies analytically on \(\varepsilon\) for small \(\varepsilon\), and has distinct eigenvalues for \(\varepsilon = 0\), so do eigenvalues and eigenvectors of this operator, and thus the result remains valid when \(\varepsilon\) is positive sufficiently small (at fixed \(\xi\)), and we have now an eigenvector \(b_0\) of the induction equation for small enough \(\varepsilon\).

This concludes the proof of theorem \ref{theoreme_lin}.

\begin{remarque}
	Most of what we said stays valid if we take \(R_m\) out of a discrete set of values instead of only ``small enough''. Only the fact that \(\alpha(U)\) is as near as we want to a matrix of distinct non-zero eigenvalues remains a problem in that case.
	\end{remarque}

Pondering on theorem \ref{theoreme_lin}, we will prove in the next section the instability of the full nonlinear (MHD) system \eqref{MHD}.

\section{Nonlinear instability}
We will now prove nonlinear instability of the system \eqref{MHD}. For that we will adapt a method developped by S. Friedlander, W. Strauss and M. Vishik in \cite{friedlander1997nonlinear} (itself adapted from Y. Guo and W. Strauss in \cite{guo1995nonlinear}) to prove nonlinear instability given linear instability. The condition \(s>3/2+1\) in the main theorem is used only in this section to obtain an inequality on the nonlinear terms similar to that in the article. 
Actually, using the parabolic feature of the MHD system, it would be possible to  lower our regularity requirements (see for instance Friedlander et al \cite{friedlander2006nonlinear} in the context of the Navier-Stokes equations). We stick here to regular data for simplicity of exposure. 

Let \(U_s\in\ms P\) belonging to the dense subset of linearly unstable flows described in the previous section.

We introduce Leray operator \(\mb P\) to get rid of \(p\). We rewrite system \eqref{MHDPert} as 
\begin{equation}
  \label{MHDPert2}
	\partial_t\ms U + L_s\ms U = Q(\ms U,\ms U),
	\end{equation}
where
\begin{equation}
	\ms U = \begin{pmatrix} u \\ b \end{pmatrix},
	\end{equation}
\begin{equation}
	L_s\ms U = \begin{pmatrix} \mb P(u\cdot\nabla U_s + U_s\cdot\nabla u) - \frac1{R_e}\Delta u \\
		-\nabla\wedge(U_s\wedge b) - \frac1{R_m}\Delta b\end{pmatrix}
\end{equation}
is the linearized operator around \(\begin{pmatrix} U_s\\0\end{pmatrix}\),
and
\begin{equation}
	Q(\ms U,\ms U) = \begin{pmatrix} \mb P(\nabla\wedge b\wedge b - u\cdot\nabla u) \\ \nabla\wedge(u\wedge b)\end{pmatrix}
	\end{equation}
is the nonlinearity.

We shall consider this system in \(\mc H_T^s\), \(T\) to be specified later. We also denote \(L^2_T\) the \(L^2\) space on the torus. Since we will work in these spaces until the end of the proof, we will drop the \(T\) in the notation. We have the following a priori estimates for the nonlinear term:

Every term of every component of \(Q(\ms U,\ms U)\) can be written as a product of one component of \(\ms U\) and one of \(\nabla\ms U\) (composed by \(\mb P\)). Thus, we have the first estimate
\begin{equation}
	\|Q(\ms U,\ms U)\|_{L^2}\le C \|\ms U\cdot\nabla\ms U\|_{L^2} \le C \|\ms U\|_{L^2}\|\nabla\ms U\|_{L^\infty}.
	\end{equation}
For \(r = \frac12(s+n/2+1) > n/2+1\) (here \(n=3\)), we have that
\begin{equation}
	\|\nabla\ms U\|_{L^\infty}\le \|\ms U\|_{\mc H^r},
	\end{equation}
and
\begin{equation}
	\|\ms U\|_{\mc H^r} \le C \|\ms U\|_{L^2}^\eta\|\ms U\|_{\mc H^s}^{1-\eta}
	\end{equation}
for \(\eta\) such that \(r = (1-\eta)s\) (i.e. \(\eta = \frac12 - \frac{n+2}{4s} = \frac12 - \frac5{4s}\)).

Thus, we have:
\begin{equation}
	\|Q(\ms U,\ms U)\|_{L^2}\le C \|\ms U\|_{L^2}^{1+\eta}\|\ms U\|_{\mc H^s}^{1-\eta}
	\end{equation}



We shall prove the following theorem, which is a mere reformulation of the first theorem:
\begin{theoreme}
	Suppose \(s>3/2+1\), \(R_m\) small enough and \(R_e>0\). Then there exists a vector \(T = (T_1,T_2,T_3)\) with integer coordinates such that the system \eqref{MHDPert2} is nonlinearly unstable in \(\mc H_T^s\) in the following sense:

	There exists a growing mode \(\ms U_0=\begin{pmatrix} u_0 \\ b_0\end{pmatrix}\) with \(\|\ms U_0\|_{\mc H_T^s} = 1\) of the linearized system and a constant \(C_0>0\) depending only on \(\ms U_0,s,n,\rho\) such that for all \(\delta>0\), the perturbation of the linear solution \(\delta\ms U_0\) satisfies \(\|\ms U(t_\delta)\|_{\mc H_T^s}\ge C_0\) for some time \(t_\delta\). Furthermore, \(t_\delta\sim_{\delta\rightarrow0} -\frac{\ln\delta}{\rho}\).
	\end{theoreme}
\begin{remarque}
	This theorem only proves that the solution reaches high values in finite time however small the initial condition is. In particular it doesn't say wether it explodes in finite time or is defined for all time, since both situations are possible in the result of the theorem. Actually, for small enough initial data (that is small enough \(\delta\)), it is classical that the solution of the system doesn't explode in finite time and thus the approximation of \(t_\delta\) for vanishing \(\delta\) is valid. More details in \cite{temam2001navier} (Theorem 3.8)
	\end{remarque}

Assume in contrary that the system is nonlinearly stable while spectrally unstable, that is that for any growing mode of initial value \(\ms U_0=\begin{pmatrix} u_0 \\ b_0\end{pmatrix}\) with \(\|\ms U_0\|_{\mc H_T^s} = 1\) in any box \(T\) of the system (i.e. an eigenvector of \(L_s\) associated to an eigenvalue with positive real part) and for any \(\varepsilon>0\), there exists a \(\delta>0\) such that for all \(t>0\), the solution \(\ms U(t)\) of the system with initial values \(\delta\ms U_0\) satisfies \(\|\ms U(t)\|_{\mc H_T^s}\le\varepsilon\).



\subsection{Proof of the nonlinear instability}

Let \(b_L(\theta) = e^{\iim\varepsilon\xi\cdot\theta}b^\varepsilon(\theta)\) be the unstable eigenmode given by Theorem \(2\). We fix \(\varepsilon<\varepsilon_0\) in such a way that 
\begin{equation}
  T = \left(\frac{2\pi}{\varepsilon|\xi_1|},\frac{2\pi}{\varepsilon|\xi_2|},\frac{2\pi}{\varepsilon|\xi_3|}\right)
\end{equation}
belongs to \(\mb N^3\). This is possible as \(\xi\in(2\pi\mb Q_+^*)^3\) (take \(\varepsilon\) to be the inverse of a large integer). In particular, \(\tilde{\ms U_0} = \begin{pmatrix} 0\\b_L\end{pmatrix}\) is an eigenvector of \(L_s\) associated to an eigenvalue with positive real part in \(\mc H_T^s\).

Denote by \(\rho\) the maximum real part of the spectrum of \(L_s\):
\begin{equation}
  \rho = \max\left\{\Re\lambda,\lambda\in\spectr(L_s)\right\}.
  \end{equation}
By the previous remark on \(\tilde{\ms U_0}\), we know that \(\rho>0\). Moreover, since the spectrum of \(L_s\) is only made of eigenvalues, there exists an eigenvector \(\ms U_0 = \begin{pmatrix} u_0\\b_0\end{pmatrix}\) with eigenvalue \(\lambda\) satisfying exactly \(\Re\lambda = \rho\).

Moreover, by standard properties of the spectral radius of \(e^{L_s}\), we know that
\begin{equation}
  e^\rho = \lim_{t\rightarrow+\infty}\|e^{tL_s}\|^{1/t}.
  \end{equation}
Thus, for all \(\eta>0\), there exists a \(C_\eta\) such that for all \(t\ge0\),
\begin{equation}
  \|e^{tL_s}\|\le C_\eta e^{\rho(1+\frac\eta2)t}
  \end{equation}

Then, the solution \(\ms U\) of \eqref{MHDPert2} with initial data \(\delta U_0\) can be written
\begin{equation}
  \ms U = \dot{\ms U}+\delta\ms U_0,
  \end{equation}
where \(\dot{\ms U}\) satisfies
\begin{equation}
  \left\{\begin{array}{l}
  \partial_t\dot{\ms U} = L_s\dot{\ms U} + Q(\ms U,\ms U)\\
  \dot{\ms U}_{|t=0}=0
  \end{array}\right.
  \end{equation}

Using Duhamel's formula, we can write the solution under the form:
\begin{equation}
	\dot{\ms U}(t) = \int_0^te^{(t-s)L_s}Q(\ms U,\ms U)\ud s.
	\end{equation}
	
Define
\begin{equation}
	T_\delta = \sup\left\{t>0,\forall s\le t,\|\dot{\ms U}(s)\|_{L^2}\le\frac12\delta e^{\rho s}\|\ms U_0\|_{L^2}\right\},
	\end{equation}
where \(\lambda\) is the coefficient in the growing linear mode. Then, since we supposed that \(\ms U = \ms U_l + \dot{\ms U}\) remains bounded, and \(\|\ms U_l\|\sim \delta e^{\rho t}\|\ms U_0\|\), we have that \(T_\delta < \infty\).

And for all \(t<T_\delta\),
\begin{equation}
	\|\ms U(t)\|_{L^2}\le \frac32\delta e^{\rho t}\|\ms U_0\|_{L^2}.
	\end{equation}

Using a priori estimates stated earlier, we have for \(\eta = \frac12 - \frac{n+2}{4s} = \frac12 - \frac{5}{4s}\):
\begin{equation}
	\|\dot{\ms U}(t)\|_{L^2} \le C_\eta \int_0^te^{\rho(1+\frac\eta2)(t-s)}\|\ms U\|_{L^2}^{1+\eta}\|\ms U\|_{\mc H^s}^{1-\eta}\ud s.
	\end{equation}
\begin{equation}
	\begin{split}
	\|\dot{\ms U}(t)\|_{L^2} &\le \delta^{1+\eta}C_\eta \int_0^te^{(t-s)(1+\frac\eta2)\rho}e^{(1+\eta)\rho s}\varepsilon^{1-\eta}\ud s \\
		&\le C_\eta\delta^{1+\eta}e^{t(1+\frac\eta2)\rho}\varepsilon^{1-\eta}\frac2{\eta\rho}(e^{t\frac\eta2\rho}-1) \\
		&\le C'\delta^{1+\eta}e^{t(1+\eta)\rho}\varepsilon^{1-\eta}\\
		\end{split}
	\end{equation}
Thus, getting back to \(\ms U\):
\begin{equation}
	\|\ms U\|_{L^2}\ge\delta\|\ms U_0\|_{L^2}e^{\rho t} - C'\delta^{1+\eta}\varepsilon^{1-\eta}e^{t(1+\eta)\rho}.
	\end{equation}
By definition of \(T_\delta\), we have:
\begin{equation}
	\|\dot{\ms U}(T_\delta)\|_{L^2} = \frac12\delta e^{\rho T_\delta}\|\ms U_0\|_{L^2}\le C'\delta^{1+\eta}\varepsilon^{1-\eta}e^{T_\delta(1+\eta)\rho},
	\end{equation}
that is
\begin{equation}
	\frac12\|\ms U_0\|_{L^2}\le C'\delta^\eta\varepsilon^{1-\eta}e^{T_\delta\eta\rho},
	\end{equation}
\begin{equation}
	T_\delta\ge\ln\left(\frac{\|\ms U_0\|_{L^2}}{2C'\delta^\eta\varepsilon^{1-\eta}}\right)\frac1{\eta\rho}.
	\end{equation}
Denote by \(t_\delta\) the quantity
\begin{equation}
	t_\delta = \ln\left(\frac{\|\ms U_0\|_{L^2}}{2C'\delta^\eta\varepsilon^{1-\eta}(1+\eta)}\right)\frac1{\eta\rho}.
	\end{equation}
Thus \(t_\delta \le T_\delta\), and at this time \(t_\delta\)
\begin{equation}
	\|\ms U\|_{L^2} \ge \delta e^{\rho t_\delta}\left(\|\ms U_0\|_{L^2} - \frac{\|\ms U_0\|_{L^2}}{1+\eta}\right) = \frac\eta{1+\eta}\|\ms U_0\|_{L^2}\left(\frac{\|\ms U_0\|_{L^2}}{C'\varepsilon^{1-\eta}(1+\eta)}\right)^{1/\eta}.
	\end{equation}
Thus \(\|\ms U\|_{L^2}\) (and \(\|\ms U\|_{\mc H^s}\)) can be bounded below by a value independent of \(\delta\), which contradicts the asumption, and proves the theorem.

\subsection{Estimation of the growing time}
We will now estimate the time in which \(\|\ms U\|_{\mc H^s}\) becomes of order \(1\). Denote by \(T_0\) this time:
\begin{equation}
	T_0 = \sup\{t>0,\forall s\le t,\|\ms U(s)\|_{\mc H^s}\le1\}.
	\end{equation}
We use again the definition of \(T_\delta\):
\begin{equation}
	T_\delta = \sup\{t>0,\forall s\le t,\|\dot{\ms U}(s)\|_{L^2}\le\frac12\delta e^{\rho s}\|\ms U_0\|_{L^2}\}.
	\end{equation}
This time, \(T_\delta\) might be infinite. For big enough \(\lambda\) as in the former proof, we have again that
\begin{equation}
	T_\delta \ge\ln\left(\frac{\|\ms U_0\|_{L^2}}{2C'\delta^\eta}\right)\frac1{\eta\rho}
	\end{equation}
Denote again by \(t_\delta\) the right-hand side:
\begin{equation}
	t_\delta = \ln\left(\frac{\|\ms U_0\|_{L^2}}{2C'\delta^\eta(1+\eta)}\right)\frac1{\eta\rho}.
	\end{equation}
Thus in time \(t_0 = \inf\{t_\delta,T_0\}\) and with the constant 
\begin{equation}
	C_0 = \inf\left\{1,\frac\eta{1+\eta}\|\ms U_0\|_{L^2}\left(\frac{\|\ms U_0\|_{L^2}}{C'(1+\eta)}\right)^{1/\eta}\right\},
	\end{equation}
we have
\begin{equation}
	\|\ms U(t_0)\|_{L^2}\ge C_0,
	\end{equation}
and 
\begin{equation}
	t_0 \le \ln\left(\frac{\|\ms U_0\|_{L^2}}{2C'(1+\eta)}\right)\frac1{\eta\rho} -\frac{\ln\delta}{\rho}\le \frac{-\ln\delta}{\rho}+C_0'.
	\end{equation}


\begin{acknowledgements}
	The author acknowledges the support of the ANR Grant 08-JCJC-0104, coordinated by David G\'erard-Varet, and thanks the latter for suggesting the subject, reading the manuscript and many advices to achieve the result.
	\end{acknowledgements}

\clearpage
\bibliographystyle{abbrv}
\bibliography{biblio_perso} 
\end{document}